\documentclass[12pt,draftcls,onecolumn]{IEEEtran}                                                      % if you need a4paper
\IEEEoverridecommandlockouts                              % This command is only
\usepackage{amsfonts}
\usepackage{graphicx}
\usepackage{amsmath}
\usepackage{amssymb}
\usepackage{geometry}
\usepackage{fancyhdr}
\usepackage{url}
\usepackage{epsfig}
\usepackage{epstopdf}

\usepackage{enumerate}

\newtheorem{theorem}{Theorem}

\newtheorem{remark}{Remark}

\title{\LARGE \bf
Separation of Quantization and Control in \\  Optimal Analog to Digital Converters\thanks{}}

\author{Mitra Osqui$\dagger $ \ \ \ \ \ \ \ Alexandre Megretski$\ddagger $ %\thanks{} 
\thanks{$\dagger$Mitra Osqui is currently a Ph.D. candidate at the department of EECS,
Laboratory for Information and Decision Systems (LIDS) at the Massachusetts
Institute of Technology, Cambridge, MA. E-mail: mitra@mit.edu}\thanks{$%
\ddagger $ Alexandre Megretski is currently a professor of EECS at LIDS at
MIT, Cambridge, MA. E-mail: ameg@mit.edu. }}

\begin{document}

\maketitle
\thispagestyle{empty}
\pagestyle{empty}

%%%%%%%%%%%%%%%%%%%%%%%%%%%%%%%%%%%%%%%%%%%%%%%%%%%%%%%%%%%%%%%%%%%%%%%%%%%%%%%%
\begin{abstract}

In this paper we prove optimality of a certain class of Analog to Digital Converters (ADCs), which can be viewed as generalized Delta-Sigma Modulators (DSMs), with respect to a performance measure that can be characterized as the worst-case average intensity of the signal representation error. An analytic expression for the ADC performance is given. Furthermore, our result proves separation of quantization and control for this class of ADCs subject to some technical conditions.

%the set of possible output values are uniformly spaced with spacing $\delta$, a simple analytic expression for its performance. 

\end{abstract}
 %
%%%%%%%%%%%%%%%%%%%%%%%%%%%%%%%%%%%%%%%%%%%%%%%%%%%%%%%%%%%%%%%%%%%%%%%%%%%%%%%%
\section{INTRODUCTION AND MOTIVATION}

Analog to Digital Converters (ADCs) act as the interface between the analog world and digital processors. 
They are present in almost all digital control and communication systems and modern high-speed data 
conversion and storage systems. Naturally, the design and analysis of ADCs have, for many years, attracted 
the attention and interest of researchers from various disciplines across academia and industry. Despite 
the progress that has been made in this field, the design of optimal ADCs remains an open challenging problem, 
and the fundamental limitations of their performance are not well understood. This paper is concerned 
with the latter problem.
\thispagestyle{empty}\pagestyle{empty}

A particular class of ADCs primarily used in high resolution applications
is\ the Delta-Sigma Modulator (DSM). Fig. \ref{fig_t} illustrates the
classical first-order DSM \cite{PhDBib:Oppenheim}, where $Q$ is a
quantizer with uniform step size.

%%%%%%%%%%%% FIGURE %%%%%%%%%%%%%%%%%
\setlength{\unitlength}{2200sp}
\begingroup\makeatletter\ifx\SetFigFont\undefined\gdef\SetFigFont#1#2#3#4#5{\reset@font\fontsize{#1}{#2pt}\fontfamily{#3}\fontseries{#4}\fontshape{#5}\selectfont}\fi\endgroup
\begin{figure}[h]\begin{center}\begin{picture}(6549,1524)(514,-3298)
\thinlines
\put(1351,-2161){\circle{336}}
\put(2326,-2536){\framebox(1125,750){}}
\put(4726,-2536){\framebox(1125,750){}}

\put(3500,-3700){\framebox(1125,750){}}

\put(1501,-2161){\vector( 1, 0){825}}
\put(3451,-2161){\vector( 1, 0){1275}}
\put(6601,-2161){\line( 0,-1){1125}}

\put(6601,-3286){\vector(-1, 0){2000}}

\put(3485,-3286){\line(-1, 0){2130}}

\put(1351,-3286){\vector( 0, 1){975}}
\put(526,-2161){\vector( 1, 0){675}}
\put(5851,-2161){\vector( 1, 0){1200}}
\put(2500,-2236){\makebox(0,0)[lb]{\smash{\SetFigFont{12}{14.4}{\rmdefault}{\mddefault}{\updefault}$\frac{1}{1-z^{-1}}$}}}

\put(3800,-3436){\makebox(0,0)[lb]{\smash{\SetFigFont{12}{14.4}{\rmdefault}{\mddefault}{\updefault}$z^{-1}$}}}

\put(1470,-2536){\makebox(0,0)[lb]{\smash{\SetFigFont{12}{14.4}{\rmdefault}{\mddefault}{\updefault}-}}}
\put(850,-2400){\makebox(0,0)[lb]{\smash{\SetFigFont{10}{14.4}{\rmdefault}{\mddefault}{\updefault}$+$}}}
%\put(1215,-2250){\makebox(0,0)[lb]{\smash{\SetFigFont{12}{14.4}{\rmdefault}{\mddefault}{\updefault}$+$}}}
\put(450,-2011){\makebox(0,0)[lb]{\smash{\SetFigFont{10}{14.4}{\rmdefault}{\mddefault}{\updefault}$r[n]$}}}
\put(3800,-2011){\makebox(0,0)[lb]{\smash{\SetFigFont{10}{14.4}{\rmdefault}{\mddefault}{\updefault}$y[n]$}}}
\put(5100,-2236){\makebox(0,0)[lb]{\smash{\SetFigFont{14}{14.4}{\rmdefault}{\mddefault}{\updefault}$Q$}}}
\put(6350,-2011){\makebox(0,0)[lb]{\smash{\SetFigFont{10}{14.4}{\rmdefault}{\mddefault}{\updefault}$u[n]$}}}

\end{picture}
\end{center}
\vspace{.05in}
\caption{Classical First-Order Sigma-Delta Modulator}
\vspace{-.05in}
\label{fig_t}
\end{figure}
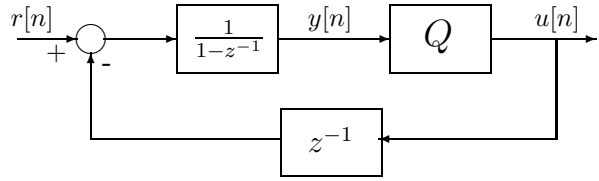
%%%%%%%%%%%%%%%%%%%%%%%%%%%%%%%%%%%%%%%%%%%

An extensive body of research on DSMs has appeared in the signal processing
literature. One well known approach is based on linearized additive noise
models and filter design for noise shaping \cite{PhDBib:Oppenheim}-\nocite%
{Derpich2008}\nocite{LNM1}\nocite{PhDBib:adhocdesign2}\nocite{PhDBib:UnderstandingDSMs}\cite{PhDBib:Norsworthy1997}. The underlying assumption for validity of the
linearized additive noise model is availability of a relatively high number
of bits. Alternative approaches based on a formalism of the signal
transformation performed\ by the quantizer have been exploited for
deterministic analysis in \cite{PhDBib:ThaoVetterli31}-\nocite%
{PhDBib:ThaoVetterli30}\cite{PhDBib:Thao2006Tile}. Some other works that do
not use linearized additive noise models are reported in \cite{QuevedoGoodwin}-\nocite%
{SteinerYang}\cite{Wang}.

In control literature, \cite{BOYDWang2}-\nocite{BOYDWang}\cite{BOYDWang3} find performance bounds and suboptimal policies for linear stochastic control problems using Bellman inequalities with quadratic value functions. The problem is relaxed and solved using linear matrix inequalities and semidefinite programming. For references on quantized control, please see \cite{Bullo2006}-\nocite{Brockett2000}\cite{EliaMitter}. 
%Work on optimal dynamic quantization can be found in \cite{Azuma2008J}.

In \cite{PhDBib:Mitra1} and \cite{PhDBib:Mitra2} we provided a characterization of the solution to the optimal ADC design problem and presented a generic methodology for
numerical computation of sub-optimal solutions along with computation of a certified upper bound and lower bound on the performance, respectively.

Fig. \ref{fig0} illustrates the setup we use for measuring the performance of the ADC. The performance of an ADC is evaluated with respect to a cost function which is a measure of the intensity of the error signal $e$ (the difference between the input signal $r$ and its quantized version $u$) for the worst case input sequence. The error signal is passed through a shaping filter which dictates the frequency region in which the error is to be minimized. Furthermore, we show that the dynamical system within the optimal ADC is a copy of the shaping filter used to define the performance criteria.

\setlength{\unitlength}{2050sp}\begingroup\makeatletter\ifx\SetFigFont
\undefined\gdef\SetFigFont#1#2#3#4#5{ \reset@font\fontsize{#1}{#2pt}
\fontfamily{#3}\fontseries{#4}\fontshape{#5} \selectfont}\fi\endgroup
\begin{figure}[h]
\begin{center}
\begin{picture}(5649,1599)(289,-2773)
\thinlines\put(3850,-1636){\circle{212}}
\put(1500,-2086){\framebox(1350,900){}}
\put(4576,-2086){\framebox(1350,900){}}
\put(2850,-1636){\vector( 1, 0){900}}
\put(3950,-1636){\vector( 1, 0){650}}
\put(5926,-1636){\vector( 1, 0){840}}
\put(3850,-2761){\vector( 0, 1){1000}}
\put(900,-2761){\line(1, 0){2950}}
\put(900,-2761){\line( 0, 1){1125}}
\put(450,-1636){\vector( 1, 0){1080}}
\put(1820,-1750){\makebox(0,0)[lb]{\smash{\SetFigFont{12}{16.8}{\rmdefault}{\mddefault}{\updefault}ADC}}}
\put(3950,-2000){\makebox(0,0)[lb]{\smash{\SetFigFont{10}{16.8}{\rmdefault}{\mddefault}{\updefault}$+$}}}
\put(3400,-1850){\makebox(0,0)[lb]{\smash{\SetFigFont{10}{16.8}{\rmdefault}{\mddefault}{\updefault}$-$}}}
\put(600,-1450){\makebox(0,0)[lb]{\smash{\SetFigFont{10}{16.8}{\rmdefault}{\mddefault}{\updefault}$r$}}}
\put(3200,-1450){\makebox(0,0)[lb]{\smash{\SetFigFont{10}{16.8}{\rmdefault}{\mddefault}{\updefault}$u$}}}
\put(4200,-1450){\makebox(0,0)[lb]{\smash{\SetFigFont{10}{16.8}{\rmdefault}{\mddefault}{\updefault}$e$}}}
\put(6500,-1450){\makebox(0,0)[lb]{\smash{\SetFigFont{10}{16.8}{\rmdefault}{\mddefault}{\updefault}$q$}}}
\put(4680,-1550){\makebox(0,0)[lb]{\smash{\SetFigFont{14}{16.8}{\rmdefault}{\mddefault}{\updefault}{\normalsize Shaping}}}}
\put(4890,-1950){\makebox(0,0)[lb]{\smash{\SetFigFont{14}{16.8}{\rmdefault}{\mddefault}{\updefault}{\normalsize Filter}}}}
\end{picture}
\end{center}
\par
\vspace*{0.0in} \vspace*{-0.1in}\caption{Setup Used for Measuring the
Performance of the ADC}%
\label{fig0}%
\end{figure}
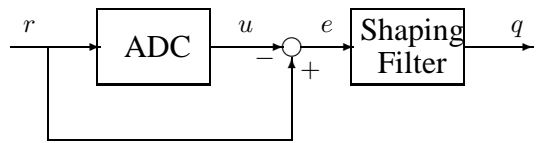

 In \cite{PhDBib:Mitra1} we also presented an exact analytical solution to the optimal ADC for first-order shaping filters, and showed that the classical first-order DSM (Figure \ref{fig_t}) is identical to our optimal ADC. This result proved the optimality of the classical first-order DSM with respect to the adopted performance measure, and was a step towards understanding the limitations of performance. In this paper we provide the optimal solution for higher order shaping filters subject to certain technical conditions and prove optimality of some higher order DSMs.

\subsection*{Notation and Terminology:}
\begin{itemize}

\item Given a set $P$, $\ell_+(P)$ is the set of all one-sided sequences $x$ with values in $P$, i.e. functions $x:\mathbb{Z}_+ \mapsto P$.
%Given a set $P$, $\ell_+(P)$ is the set of all sequences that map $\mathbb{Z}_+$ to $P$:
%\begin{align}
%\ell_+(P) & \;  {\buildrel\rm def\over =}\; \left\{x:\mathbb{Z}_+ \mapsto P \right\}. 
%\end{align}

%\medskip
%
%\item The $\infty-$norm is defined as:
%\begin{equation*}
%\|v\|_\infty = \max |v_i|, \quad \text{for} \quad v =  \left[
%\begin{array}
%[c]{ccc}
%v_1 \\  \vdots \\  v_m
%\end{array}
%\right] \in  \mathbb{R}^m
%\end{equation*}
%and
%\begin{equation*}
%\|M\|_\infty  \;  {\buildrel\rm def\over =}\; \sup_{v \ne 0}\frac{\|Mv\|_\infty}{\|v\|_\infty} =\max_{i \in\{1, \cdots, l\}} \sum_{j=1}^{m}|M_{ij}|
%\end{equation*}
%for a matrix $M = \left(M_{ij}\right) \in \mathbb{R}^{l \times m}$.

\end{itemize}
 %
%%%%%%%%%%%%%%%%%%%%%%%%%%%%%%%%%%%%%%%%%%%%%%%%%%%%%%%%%%%%%%%%%%%%%%%%%%%%%%%%

\section{Problem Formulation\label{ProbForm}}

The problem setup in this section is taken from \cite{PhDBib:Mitra1}.

%%%%%%%%%%%%
\subsection{Analog to Digital Converters}

In this paper, a general ADC is viewed as a causal, discrete-time, non-linear
system $\Psi,$ accepting arbitrary inputs in the $[-1,1]$ range and producing
outputs in a fixed finite subset $U\subset\mathbb{R},$ as shown in Fig. \ref{fig_a}. We assume $\max U > 1$ and $\min U < -1$.\vspace{0.15in}%

\setlength{\unitlength}{1800sp}
\begingroup\makeatletter\ifx\SetFigFont\undefined\gdef\SetFigFont
#1#2#3#4#5{\reset@font\fontsize{#1}{#2pt}\fontfamily{#3}\fontseries
{#4}\fontshape{#5}\selectfont}\fi\endgroup\begin{figure}[h]\begin{center}%
\begin{picture}(6549,1524)(514,-3298)
\thinlines\put(2350,-2736){\framebox(2200,1250){}}
\put(800,-2161){\vector( 1, 0){1570}}
\put(4550,-2161){\vector( 1, 0){1275}}
\put(3200,-2300){\makebox(0,0)[lb]{\smash{\SetFigFont{18}{14.4}{\rmdefault}{\mddefault}{\updefault}$\Psi$}}}
\put(200,-1900){\makebox(0,0)[lb]{\smash{\SetFigFont{10}{14.4}{\rmdefault}{\mddefault}{\updefault}$r[n] \in[-1,1]$}}}
\put(200,-2560){\makebox(0,0)[lb]{\smash{\SetFigFont{10}{14.4}{\rmdefault}{\mddefault}{\updefault}$n \in\mathbb{Z}_+$}}}
\put(5000,-1900){\makebox(0,0)[lb]{\smash{\SetFigFont{10}{14.4}{\rmdefault}{\mddefault}{\updefault}$u[n]\in U$}}}
\put(5000,-2560){\makebox(0,0)[lb]{\smash{\SetFigFont{10}{14.4}{\rmdefault}{\mddefault}{\updefault}$n \in\mathbb{Z}_{+}$}}}
\end{picture}
\end{center}
\vspace*{-.2in}
\caption{Analog to Digital Converter as a Dynamical System}
\vspace{-.05in}
\label{fig_a}
\end{figure}
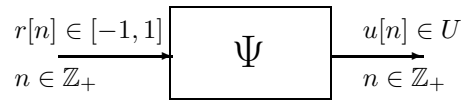%

Equivalently, an ADC is defined by a sequence of functions
$\Upsilon_{n}:\left[  -1,1\right]  ^{n+1}\mapsto U$ according to
\smallskip
\begin{equation}
\Psi:~u[n]=\Upsilon_{n}\left(  r[n]  ,r[ n-1], \cdots ,r[0]  \right)
,~~n\in\mathbb{Z}_{{\small +}}.\label{general_u}%
\end{equation}

The class of ADCs defined above is denoted by $\mathcal{Y}_{U}.$

%%%%%%%%%%%%%%%%%%%%%%
\subsection{Asymptotic Weighted Average Intensity (AWAI) of a Signal}

Let $\phi:\mathbb{R}\mapsto\mathbb{R}_{+}$ be an even, non-negative, and monotonically nondecreasing function on the positive reals; and $G\left(  z\right)$ be the transfer function of a strictly causal LTI dynamical system $L_{G}$ with input $w$ and output $q$: 
\begin{equation}
L_G:
\begin{cases}
x[n+1]  =Ax[n]  +Bw[n], \quad x[0] =0, \label{LG0}\\ q[n]  =Cx[n]
\end{cases}
\end{equation}

\noindent where $A,$ $B,$ $C$ are given matrices of appropriate dimensions. The Asymptotic Weighted Average Intensity $\eta_{G,\phi}\left(  w\right)$ of signal $w$ with respect to $G\left(  z\right)$  and $\phi$  is given by:
\begin{equation}
\eta_{G,\phi}\left(  w\right)  =\underset{N\mapsto\infty}{\lim\sup}%
\frac{1}{N}
{\displaystyle\sum\limits_{n=0}^{N-1}}
\phi\left(  q[n]  \right).\label{AWA}%
\end{equation}
 Examples of functions $\phi$ to consider are: $\phi(q )=\left\vert q \right \vert$ and $\phi(q)=\left\vert q \right \vert^2$. We assume without loss of generality that $CB \ne 0$. Indeed, since $\eta_{G,\phi}$ does not change if $G(z)$ is replaced by $zG(z)$, i.e. if $q[n]$ is replaced with $q[n+1]$ in \eqref{LG0}, the case when $CB=0$ can be reduced to the case $CB \ne 0$ by extracting a delay from $L_G$.

 %

%\setlength{\unitlength}{2200sp}\begingroup\makeatletter\ifx\SetFigFont
%\undefined\gdef\SetFigFont#1#2#3#4#5{  \reset@font\fontsize{#1}{#2pt}
%\fontfamily{#3}\fontseries{#4}\fontshape{#5}  \selectfont}\fi\endgroup
%\begin{figure}[h]\begin{center}\begin{picture}(8649,1599)(289,-2773)
%\thinlines\put(2576,-2086){\framebox(1350,900){}}
%\put(1550,-1636){\vector( 1, 0){1000}}
%\put(3926,-1636){\vector( 1, 0){1000}}
%\put(1600,-1486){\makebox(0,0)[lb]{\smash{\SetFigFont{10}{16.8}{\rmdefault}{\mddefault}{\updefault}$w[n]$}}}
%\put(4250,-1486){\makebox(0,0)[lb]{\smash{\SetFigFont{10}{16.8}{\rmdefault}{\mddefault}{\updefault}$q[n]$}}}
%\put(3000,-1700){\makebox(0,0)[lb]{\smash{\SetFigFont{14}{16.8}{\rmdefault}{\mddefault}{\updefault}$L_G$}}}
%\end{picture}
%\end{center}
%\vspace*{-.3in}
%\caption{Strictly Proper LTI Shaping Filter $G(z)$}
%\vspace*{-.2in}
%\label{fig_b}
%\end{figure}
%%%%%%%%%%%%%%
  %
\subsection{ADC Performance Measure\label{PerfMeas}}

The setup that we use to measure the performance of an ADC is illustrated in
Fig. \ref{fig1}. The performance measure of $\Psi\in\mathcal{Y}_{U}$,
denoted by $\mathcal{J}_{G,\phi}\left(  \Psi\right)  ,$ is the worst-case AWAI
of the error signal for all input sequences $r \in \ell_+(\left[
-1,1\right])  ,$ that is:%
\begin{equation}
\mathcal{J}_{G,\phi}\left(  \Psi\right)  =\sup_{r\in\ell_+(\left[-1,1\right])}\eta_{G,\phi
}\left( r - \Psi\left(  r\right)\right). \label{perfmeas}
\end{equation}

\setlength{\unitlength}{2200sp}\begingroup\makeatletter\ifx\SetFigFont
\undefined\gdef\SetFigFont#1#2#3#4#5{  \reset@font\fontsize{#1}{#2pt}
\fontfamily{#3}\fontseries{#4}\fontshape{#5}  \selectfont}\fi\endgroup
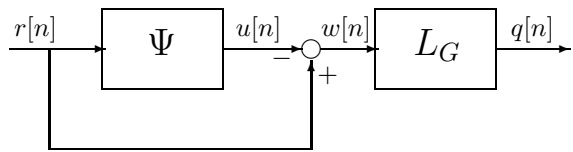
\begin{figure}[h]\begin{center}\begin{picture}(8649,1599)(289,-2773)
\thinlines\put(3850,-1636){\circle{212}}
\put(1500,-2086){\framebox(1350,900){}}
\put(4576,-2086){\framebox(1350,900){}}
\put(2850,-1636){\vector( 1, 0){900}}
\put(3950,-1636){\vector( 1, 0){650}}
\put(5926,-1636){\vector( 1, 0){840}}
\put(3850,-2761){\vector( 0, 1){1000}}
\put(900,-2761){\line(1, 0){2950}}
\put(900,-2761){\line( 0, 1){1125}}
\put(450,-1636){\vector( 1, 0){1080}}
\put(2000,-1700){\makebox(0,0)[lb]{\smash{\SetFigFont{14}{16.8}{\rmdefault}{\mddefault}{\updefault}$\Psi$}}}
\put(3920,-2011){\makebox(0,0)[lb]{\smash{\SetFigFont{10}{16.8}{\rmdefault}{\mddefault}{\updefault}$+$}}}
%\put(3750,-1711){\makebox(0,0)[lb]{\smash{\SetFigFont{10}{16.8}{\rmdefault}{\mddefault}{\updefault}$+$}}}
\put(3400,-1820){\makebox(0,0)[lb]{\smash{\SetFigFont{10}{16.8}{\rmdefault}{\mddefault}{\updefault}$-$}}}
\put(500,-1486){\makebox(0,0)[lb]{\smash{\SetFigFont{10}{16.8}{\rmdefault}{\mddefault}{\updefault}$r[n]$}}}
\put(3000,-1486){\makebox(0,0)[lb]{\smash{\SetFigFont{10}{16.8}{\rmdefault}{\mddefault}{\updefault}$u[n]$}}}
\put(3960,-1486){\makebox(0,0)[lb]{\smash{\SetFigFont{10}{16.8}{\rmdefault}{\mddefault}{\updefault}$w[n]$}}}
\put(6100,-1486){\makebox(0,0)[lb]{\smash{\SetFigFont{10}{16.8}{\rmdefault}{\mddefault}{\updefault}$q[n]$}}}
\put(5000,-1700){\makebox(0,0)[lb]{\smash{\SetFigFont{14}{16.8}{\rmdefault}{\mddefault}{\updefault}$L_G$}}}
\end{picture}
\end{center}
\vspace*{0.0in}
\caption{Setup Used for Measuring the Performance of the ADC}
\vspace*{-.1in}
\label{fig1}
\end{figure}

%that when $\phi\left(  \cdot\right)  =\left\vert
%\cdot\right\vert ^{2}$ and $L_{G}$ is a strictly stable dynamical system with
%transfer function $G\left(  z\right)  ,$ the AWA can be interpreted as the
%average power of the filtered input for signals which are sums of sinusoids:%
%\begin{equation}
%w\left[  n\right]  =%
%%TCIMACRO{\dsum \limits_{k=0}^{\infty}}%
%%BeginExpansion
%{\displaystyle\sum\limits_{k=0}^{\infty}}
%%EndExpansion
%w_{k}e^{j\omega_{k}n}\mapsto\eta_{G,\phi}\left(  w\right)  =%
%%TCIMACRO{\dsum \limits_{k=0}^{\infty}}%
%%BeginExpansion
%{\displaystyle\sum\limits_{k=0}^{\infty}}
%%EndExpansion
%\left\vert w_{k}\right\vert ^{2}\left\vert G\left(  e^{j\omega_{k}}\right)
%\right\vert ^{2}.\label{sumsin}%
%\end{equation}
%Therefore, the AWA allows for penalizing the input of the filter over the
%frequency range of interest via the pass-band of $G\left(  z\right)  .$ An
%alternative measure can be obtained with $\phi\left(  \cdot\right)
%=\left\vert \cdot\right\vert ,$ which is attractive due to its simplifying
%properties. In this case, the AWA represents the average amplitude of the
%filtered input signal.

%%%%%%%%%%%%%%%%
\subsection{ADC Optimization}

Given $L_{G}$ and $\phi,$ we consider $\Psi_{o}\in\mathcal{Y}_{U}$ an optimal
ADC if $\mathcal{J}_{G,\phi}\left(  \Psi_{o}\right)  \leq\mathcal{J}_{G,\phi
}\left(  \Psi\right)  $ for all $\Psi\in\mathcal{Y}_{U}.$ The corresponding
optimal performance measure $\gamma_{G,\phi}\left(  U\right)  $ is defined as%
\begin{equation}
\gamma_{G,\phi}\left(  U\right)  =\underset{\Psi\in\mathcal{Y}_{U}}{\inf
}\mathcal{J}_{G,\phi}\left(  \Psi\right).\label{cutefuzzfuzz}%
\end{equation}

%The objective is to find  $\gamma_{G,\phi}\left(  U\right) $.

%%%%%%%%%%%%%%%%%%%%%%%%%%%%%%%%%%%%%%%%%%%%%%%%%%%%%%%%%%

\section{Our Approach\label{OurApp}}

We search for the optimal ADC within the class of time invariant state-space models and associate the optimal ADC
design problem with a full-information feedback control problem. We show for a certain class of ADCs that the setup depicted in Figure \ref{figStateFeedbackUB} is an optimal ADC architecture. The function $K:\mathbb{R}^{m}\mathbb{\times }\left[  -1,1\right]  \mapsto U$ is said to be an admissible controller if there exists $\gamma \in [0,\infty)$ such that every triplet of sequences $(x_{\Psi},u,r)$ satisfying
\begin{align}
%x_{\Psi}\left[  0\right]   &  =0,\label{icxUB}\\
x_{\Psi}\left[  n+1\right]   &  =Ax_{\Psi}\left[  n\right]  +Br\left[n\right]  -Bu\left[  n\right]  ,~ x_{\Psi}\left[  0\right]   =0,\label{x+}\\
u\left[  n\right]   &  =K\left(  x_{\Psi}\left[  n\right]  ,r\left[  n\right]\right)  ,\label{OptUK}\\
q_{\Psi}\left[  n\right]   &  =Cx_{\Psi}\left[  n\right]  , \label{qPsi}%
\end{align}
\noindent also satisfies the dissipation inequality
\begin{equation}
\sup_{N,{r\in\ell_+(\left[-1,1\right])}}{\displaystyle\sum\limits_{n=0}^{N-1}}\left(  \phi\left(  q_\Psi\left[  n\right]  \right)  -\gamma\right)<\infty\label{optUconditionUB} 
\end{equation}

Note that if (\ref{optUconditionUB}) holds subject to \eqref{x+}-\eqref{qPsi}, then $\mathcal{J}_{G,\phi}\left(  \Psi\right)\leq\gamma.$ Let $\gamma_{o}$ be the maximal lower bound of $\gamma$, for which an admissible controller exists. Then $K$ is said to be an optimal controller if (\ref{optUconditionUB}) is satisfied with $\gamma=\gamma_{o}.$

\setlength{\unitlength}{2400sp}\begingroup\makeatletter\ifx\SetFigFont
\undefined\gdef\SetFigFont#1#2#3#4#5{ \reset@font\fontsize{#1}{#2pt}
\fontfamily{#3}\fontseries{#4}\fontshape{#5} \selectfont}\fi\endgroup
\begin{figure}[h]
\begin{center}
\begin{picture}(8649,2900)(289,-2773)
\thicklines\put(800,-3436){\dashbox{120}(6250,2950){}}
\thinlines\put(6470,-3100){\makebox(0,0)[lb]{\smash{\SetFigFont{10}%
{14.4}{\rmdefault}{\mddefault}{\updefault}$u[n]$}}}
\put(3201,-2850){\oval(150,150)[br]}
\put(3201,-2860){\oval(150,150)[tr]}
\put(3201,-1640){\line( 0,-1){1130}}
\put(3201,-2920){\line( 0,-1){360}}
\put(3201,-3270){\vector( 1, 0){4100}}
\put(4000,-300){\makebox(0,0)[lb]{\smash{\SetFigFont{14}{14.4}{\rmdefault
}{\mddefault}{\updefault}$\Psi$}}}
\put(930,-1000){\makebox(0,0)[lb]{\smash{\SetFigFont{10}{14.4}{\rmdefault
}{\mddefault}{\updefault}$x_{\Psi}[n]$}}}
\put(870,-1300){\vector( 1, 0){610}}
\put(870,-1300){\line( 0, 1){600}}
\put(870, -700){\line( 1, 0){5370}}
\put(6250, -700){\line( 0,-1){700}}
\put(5926,-1400){\line( 1, 0){320}}
\put(6300,-1400){\makebox(0,0)[lb]{\smash{\SetFigFont{10}{16.8}{\rmdefault
}{\mddefault}{\updefault}$x_{\Psi}[n]$}}}
\put(3850,-1636){\circle{212}}
\put(1500,-2086){\framebox(1350,900){}}
\put(4576,-2086){\framebox(1350,900){}}
\put(2850,-1636){\vector( 1, 0){900}}
\put(3950,-1636){\vector( 1, 0){635}}
\put(5926,-1800){\vector( 1, 0){775}}
\put(3850,-2855){\vector( 0, 1){1125}}
\put(900,-2855){\line(1, 0){2950}}
\put(900,-2855){\line( 0, 1){1125}}
\put(440,-1736){\vector( 1, 0){1050}}
\put(1720,-1700){\makebox(0,0)[lb]{\smash{\SetFigFont{14}{16.8}{\rmdefault}{\mddefault}{\updefault}$K(\cdot,\cdot)$}}}
\put(3650,-2050){\makebox(0,0)[lb]{\smash{\SetFigFont{10}{16.8}{\rmdefault}{\mddefault}{\updefault}$+$}}}
\put(3400,-1800){\makebox(0,0)[lb]{\smash{\SetFigFont{10}{16.8}{\rmdefault}{\mddefault}{\updefault}$-$}}}
\put(250,-1600){\makebox(0,0)[lb]{\smash{\SetFigFont{10}{16.8}{\rmdefault}{\mddefault}{\updefault}$r[n]$}}}
\put(3000,-1486){\makebox(0,0)[lb]{\smash{\SetFigFont{10}{16.8}{\rmdefault}{\mddefault}{\updefault}$u[n]$}}}
\put(3980,-1486){\makebox(0,0)[lb]{\smash{\SetFigFont{10}{16.8}{\rmdefault}{\mddefault}{\updefault}$w[n]$}}}
\put(6100,-2100){\makebox(0,0)[lb]{\smash{\SetFigFont{10}{16.8}{\rmdefault}{\mddefault}{\updefault}$q_{\Psi}[n]$}}}
\put(5000,-1700){\makebox(0,0)[lb]{\smash{\SetFigFont{14}{16.8}{\rmdefault}{\mddefault}{\updefault}$L_G$}}}
\end{picture}
\end{center}
\par
\vspace*{0.25in} \vspace*{-.0in}\caption{Full State-Feedback Control Setup}%
\label{figStateFeedbackUB}%
\end{figure}
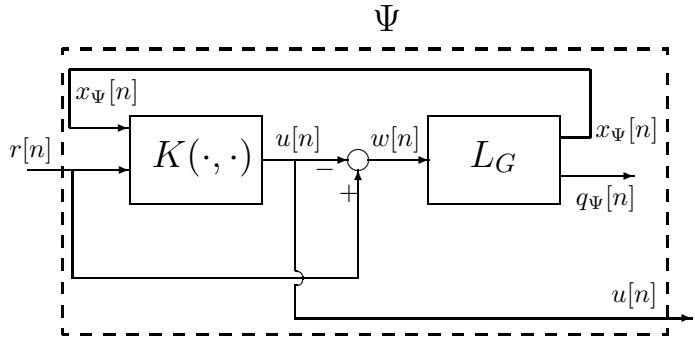

%%%%%%%%%%%%%%%%%%%%%%%%%%%%%%%%%%%%%%%%%%%
   %
\section{Main Result\label{TexSta}}

Consider the ADC optimization problem presented in Section \ref{ProbForm} with $L_G$ defined by \eqref{LG0} with $CB \ne 0$. For $\delta\in (0,2]  $ and $M\in\mathbb{N} \cup \{\infty\}$, define the set $U_M$ and function $K_M:\mathbb{R}\rightarrow U_M$ as
\begin{align}
U_M  &  =\left\{  m\delta\mathbb{~}|~m\in\mathbb{Z},~\left\vert m\right\vert
\leq M\right\} \label{SetU}\\
K_M(  \theta)   &  =\min \left \{ \arg\min_{u\in U_M}\left\vert \theta-u\right\vert \right\}.  \label{K_func}%
\end{align}

Consider the ADC $\widehat \Psi\in\mathcal{Y}_{U_M}$ defined by 
\begin{equation}
L_{\widehat \Psi} :
\begin{cases} 
x_{\widehat \Psi}\left[  n+1\right]   =Ax_{\widehat \Psi}\left[  n\right]  +Br\left[n\right]  -Bu\left[  n\right],\label{xPsi++}\\
q_{\widehat \Psi}[n]   = Cx_{\widehat \Psi}[n]\\
x_{\widehat \Psi}\left[  0\right]   =0
\end{cases}
\end{equation}
with the control law
\begin{equation}
u[n]  =K_M\left( (CB)^{-1} CAx_{\widehat \Psi}[n]  +r[n] \right). \label{Opt_u}
\end{equation}

We show in Theorem \ref{THMOptimalADC} below that if $M$ is large enough and $\delta$ is small enough, then the ADC defined above is optimal. The control decision $u[n]$ in \eqref{Opt_u} minimizes $|q_{\widehat \Psi}[n+1]|$. An interpretation of Theorem \ref{THMOptimalADC} is that a greedy algorithm is optimal subject to certain conditions. Let
\begin{equation}
q_{\widehat \Psi}[n+1] = \sum_{i=0}^k a_i q_{\widehat \Psi}[n-i] +\sum_{j=0}^k b_j (r[n-j]-u[n-j]). \label{qdiffeqn}
\end{equation}
%where $b_0=1$, which is equivalent to condition $CB=1$.
be the difference equation which is equivalent to \eqref{xPsi++}. Let $\mathcal F$ be the causal LTI system with transfer function
\begin{equation}
F(z) = \frac{1}{\displaystyle \sum_{j=0}^k b_j z^{-j}}. \label{Hz}
\end{equation}
Let $\{c_l \}_{l=0}^{\infty}$ be the unit sample response of system \eqref{qdiffeqn}, i.e.
\begin{equation}
F(z) = \displaystyle \sum_{l=0}^\infty c_l z^{-l}, \quad \text{for} \quad |z|> R_0 \label{cl}
\end{equation}
where $R_0 \in \mathbb R$ is the maximal absolute value of the largest pole of $F(z)$ in \eqref{Hz}.

\bigskip
\begin{theorem}
\label{THMOptimalADC}Let $\widehat{\Psi}\in\mathcal{Y}_{U_M}$ be the ADC defined by \eqref{xPsi++}$-$\eqref{Opt_u} with $CB \ne 0$ 
and $K_M$ defined by (\ref{SetU})$-$(\ref{K_func}). Let
\[
\beta = \left [|CB|\displaystyle \frac{\delta}{2} \left ( \displaystyle  \sum_{i=0}^k|a_i|+1 \right)+ \displaystyle \sum_{j=0}^k |b_j|\right ]\displaystyle \sum_{l=0}^\infty |c_l|,
\]
where $\{a_i \}_{i=0}^{k}$ and $\{b_j \}_{j=0}^{k}$ are defined by \eqref{qdiffeqn} and $\{c_l \}_{l=0}^{\infty}$ is defined by \eqref{Hz}$-$\eqref{cl}. Let $M \delta$ be such that $M \delta > 1$ and 
\begin{equation}
M \delta >\beta-\delta.\label{umax}
\end{equation}
 Let $f:[0,\infty)\rightarrow\lbrack 0,\infty)$ be a monotonically nondecreasing function and $\phi\left(q\right)  =f\left(  \left\vert q\right\vert \right).$ Then $\widehat{\Psi}$ is an optimal ADC in the sense that%
\begin{equation}
\mathcal{J}_{G,\phi}\left(  \Psi\right)  \geq\mathcal{J}_{G,\phi} (\widehat{\Psi})=\phi\left( |CB| \delta/2\right)  \quad\forall\Psi\in \mathcal{Y}_{U_M}.\label{PerfJ}%
\end{equation}

\end{theorem}

\bigskip

\begin{proof}
Please see the Appendix.
\end{proof}

\begin{remark}
We showed in \cite{PhDBib:Mitra1} that the first-order DSM in Figure \ref{fig_t} is optimal with respect to the shaping filter $ L_G=1/(z-1)$ with any uniform quantizer $Q$ with $M\delta > 1$.
\end{remark}
\begin{remark}
For $ L_G=z/(z-1)^2$ with any uniform quantizer $Q$ with step size $\delta \le 2$ and the magnitude of the largest value of the quantizer being larger than $1+\delta$, the second-order DSM is optimal. %The $1$-bit second-order DSM cannot satisfy this condition, hence it is not optimal. 
\end{remark}

%\begin{remark}
%An important consequence of theorem \ref{THMOptimalADC} is that a sufficient condition is that:
%if the shaping filter $L_G$ has strictly stable zeros, i.e. if $L_G$ has zeros located at $|z| \ge 1$ then $M\delta \to \infty$ for optimal design in the way i have solved this problem.
%\end{remark}
\bigskip

The optimal ADC architecture presented in Figure \ref{figStateFeedbackUB} along with the optimal control law given in \eqref{Opt_u} can be equivalently represented by Figure \ref{OptimalADC} and equation \eqref{OptimalH}, where $Q$ is a uniform quantizer with step size $\delta$ and saturation level $M\delta$ satisfying \eqref{umax} and $G(z)$ is the transfer function of the shaping filter $L_G$. Furthermore, Figure \ref{OptimalADC} has a DSM architecture, thus with a proper selection of $L_G$ as the shaping filter, many standard DSMs that satisfy the conditions in Theorem \ref{THMOptimalADC} are proven optimal. 
%\begin{figure}[h]
%\centering
%\includegraphics[scale=.9]{Figure/MyOptimalADC}
% \caption{Optimal ADC Architecture, where $G(z)=C(zI-A)^{-1}B$ is the transfer function of $L_G$}
%\label{OptimalADC}
%\end{figure}
% 
\vspace{.5in}

\setlength{\unitlength}{2500sp}\begingroup\makeatletter\ifx\SetFigFont
\undefined\gdef\SetFigFont#1#2#3#4#5{  \reset@font\fontsize{#1}{#2pt}
\fontfamily{#3}\fontseries{#4}\fontshape{#5}  \selectfont}\fi\endgroup
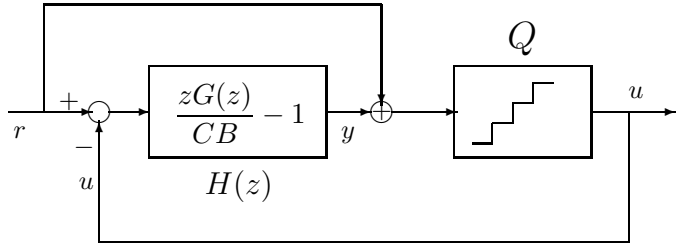
\begin{figure}[h]\begin{center}\begin{picture}(7649,1599)(289,-2773)
\thinlines\put(3850,-1636){\circle{212}}
\put(1550,-2086){\framebox(1750,900){}}
\put(4576,-2086){\framebox(1350,900){}}
\put(3300,-1636){\vector( 1, 0){450}}
\put(3950,-1636){\vector( 1, 0){650}}
\put(5926,-1636){\vector( 1, 0){840}}

\put(3850,-570){\vector( 0, -1){1000}}
\put(500,-570){\line(1, 0){3350}}
\put(500,-570){\line( 0, -1){1070}}

\put(150,-1636){\vector( 1, 0){830}}

\thinlines\put(1050,-1636){\circle{212}}

\put(1150,-1636){\vector( 1, 0){400}}

\put(1800,-1750){\makebox(0,0)[lb]{\smash{\SetFigFont{11}{16.8}{\rmdefault}{\mddefault}{\updefault}$\displaystyle\frac{zG(z)}{CB}-1$}}}
\put(2100,-2450){\makebox(0,0)[lb]{\smash{\SetFigFont{12}{16.8}{\rmdefault}{\mddefault}{\updefault}$H(z)$}}}

\put(3750,-1700){\makebox(0,0)[lb]{\smash{\SetFigFont{10}{16.8}{\rmdefault}{\mddefault}{\updefault}$+$}}}

\put(200,-1900){\makebox(0,0)[lb]{\smash{\SetFigFont{10}{16.8}{\rmdefault}{\mddefault}{\updefault}$r$}}}
\put(850,-2400){\makebox(0,0)[lb]{\smash{\SetFigFont{10}{16.8}{\rmdefault}{\mddefault}{\updefault}$u$}}}
\put(800,-2050){\makebox(0,0)[lb]{\smash{\SetFigFont{10}{16.8}{\rmdefault}{\mddefault}{\updefault}$-$}}}
\put(650,-1600){\makebox(0,0)[lb]{\smash{\SetFigFont{10}{16.8}{\rmdefault}{\mddefault}{\updefault}$+$}}}

\put(3450,-1900){\makebox(0,0)[lb]{\smash{\SetFigFont{10}{16.8}{\rmdefault}{\mddefault}{\updefault}$y$}}}

\put(6300,-1500){\makebox(0,0)[lb]{\smash{\SetFigFont{10}{16.8}{\rmdefault}{\mddefault}{\updefault}$u$}}}

\put(4750,-1950){\line( 1,0){200}}
\put(4950,-1750){\line( 0,-1){200}}

\put(4950,-1750){\line( 1,0){200}}
\put(5150,-1550){\line(0,-1){200}}

\put(5150,-1550){\line(1,0){200}}
\put(5350,-1350){\line(0,-1){200}}

\put(5350,-1350){\line(1,0){200}}

\put(5100,-1000){\makebox(0,0)[lb]{\smash{\SetFigFont{14}{16.8}{\rmdefault}{\mddefault}{\updefault}$Q$}}}

\put(6300,-1630){\line(0,-1){1300}}
\put(6300,-2930){\line( -1,0){5250}}
\put(1050,-2930){\vector(0,1){1200}}

\end{picture}
\end{center}
\caption{Optimal ADC Architecture, where $G(z)=C(zI-A)^{-1}B$ is the transfer function of $L_G$}
\label{OptimalADC}
\end{figure}

\begin{equation}
H(z) = (CB)^{-1}zG(z)-1 = (CB)^{-1}C(zI-A)^{-1}AB  \label{OptimalH}
\end{equation}

That is, if the magnitude of the largest value of the quantizer output is large enough and quantization step size is small enough, then the greedy algorithm is the optimal output for the ADC. This shows separation of quantization and control for this problem, subject to inequality \eqref{umax}.

%Higher order DSMs are not always optimal. For $ L_G=z/(z-1)^2$ with any uniform quantizer $Q$ with step size $\delta \le 2$ and the magnitude of the largest value of the quantizer being larger than $1+\delta$, the second-order DSM is optimal. The $1$-bit second-order DSM cannot satisfy this condition, hence it is not optimal. 
% 

%%%%%%%%%%%%%%%%%%%%%%%%%%%%%%%%%%%%%%%%%%%%%%%%%%%%%%%%%%%%%%%%%%%%%%%%%%%%%%%%

%\section{EXAMPLES\label{NumEx}}
%
%\subsection{Optimality of the First-Order DSM}%
%
%
%\subsection{Optimality of Higher-Order$\Delta \Sigma M$s}

%%%%%%%%%%%%%%%%%%%%%%%%%%%%%%%%%%%%%%%%%%%%%%%%%%%%%%%%%%%%%%%%%%%%%%%%%%%%%%%%

\section{CONCLUSION\label{Future}}
In this paper, we showed optimality of a certain class of ADCs (which were shown to have DSM like architecture) subject to some conditions and provided an analytic expression for the performance. We showed that there is separation of quantization and control, i.e. in the absence of quantization, the obvious choice for the optimal control law is proven to be the optimal control law given quantization, when certain technical conditions are met.

%we studied performance limitations of Analog to Digital Converters (ADCs). The performance of an ADC was defined in terms of a measure that represents the worst case average intensity of the filtered input matching error. The passband of the shaping filter defines the frequency region in which the error is to be minimized. The problem of finding a lower bound for the performance of an ADC was associated with a full information feedback optimal control problem and formulated as a dynamic game in which the input of the ADC (control variable) played against the output of the ADC (quantized disturbance). Since the disturbances can exceed the control variable in magnitude, if the shaping filter has a pole on the unit circle, then there does not exist a bounded control invariant set, which presents a challenge for numerical computations. This challenge is overcome with theoretical results that show that the value function is zero beyond a bounded region, thus computations need to be done only over this bounded region. A numerical algorithm was presented that provided certified solutions to the underlying Bellman inequality in parallel with the control law; hence, certified lower bounds on the performance of arbitrary ADCs with respect to the adopted performance criteria. 

%%%%%%%%%%%%%%%%%%%%%%%%%%%%%%%%%%%%%%%%%%%%%%%%%%%%%%%%%%%%%%%%%%%%%%%%%%%%
   %
\section{APPENDIX\label{APP}}

\subsubsection*{\bf Proof of Theorem \ref{THMOptimalADC}}
Let us begin by showing that with the control law given in \eqref{Opt_u} with $M=\infty$ we have:
\begin{equation}
\left\vert q_{\widehat{\Psi}}\left[  n\right]  \right\vert \leq|CB|\delta/2,\text{\qquad}\forall n\in\mathbb{Z}_{+}, \label{eq509}%
\end{equation}
Indeed, for $n=0,$ inequality $($\ref{eq509}$)$ follows from the initial condition in (\ref{xPsi++}). For $n > 0$,
\begin{equation*}
q_{\widehat{\Psi}}[ n+1] = CB(w[n]-K(w[n])),  \label{eq510}
\end{equation*}
where $w[n] = (CB)^{-1}CAx_{\widehat{\Psi}}[n]+r[n]$. Since $|\theta - K(\theta)| \le \delta/2$ for all $\theta \in \mathbb R$, we have \eqref{eq509} for all  $n \ge 0$. 

The next step is to use the bound $|q_{\widehat \Psi}[n]|\le |CB|\delta/2$ to show that $|u[n]| \le \beta$. Rearranging \eqref{qdiffeqn}, taking absolute value from both sides, and using the triangle inequality yields:
\begin{equation*}
\left | \displaystyle \sum_{j=0}^k b_j u[n-j] \right | \le |CB|\displaystyle \frac{\delta}{2} \left ( \displaystyle  \sum_{i=0}^k|a_i|+1 \right)+ \displaystyle \sum_{j=0}^k |b_j|\end{equation*}

If $\sum_{j=0}^k b_j u[n-j]$ is the input signal to the system $\mathcal F$ with transfer function $F(z)$ defined in \eqref{Hz}, then the output $u[n]$ is bounded in magnitude by
\begin{equation}
|u[n]| \le \beta \label{u_bound}
\end{equation}
  %
%Let $M\delta = \max \{\mu \in U_\infty : \mu < \beta \}$, i.e. $M \delta$ is the largest element of the set $U_\infty$ that is less than $\beta$. 

A sufficient condition for $|u[n]| \le M \delta$, is given by \eqref{umax}, \eqref{u_bound}, and $u \in U_\infty$. Therefore \eqref{umax} implies \eqref{eq509}.

Since both systems
$L_{G}$ and $L_{\widehat{\Psi}}$ have the same input and $x_{\widehat{\Psi}%
}\left[  0\right]  =x\left[  0\right]  =0,$ condition $(\ref{eq509})$ implies
that
\[
\left\vert q\left[  n\right]  \right\vert \leq |CB|\delta
/2,\qquad\forall n\in\mathbb{Z}_{+}.
\]
Therefore,
\[
\sup_{N,r\in[-1,1]}{\displaystyle\sum\limits_{n=0}^{N}}\left(  \phi\left( q \left[  n\right]  \right)  -\phi (|CB|\delta/2)\right)  \leq0<\infty,
\]
which implies that
\begin{equation}
\mathcal{J}_{G,\phi}(\widehat{\Psi})\leq\phi\left( |CB| \delta/2\right)  .\label{gamma_minUB}
\end{equation}
In order to complete the proof, we need to show that no ADC can achieve a
better performance than $\phi\left(  |CB|\delta/2\right)  $. It is sufficient to
show that for all $\Psi\in\mathcal{Y}_{U},$ there exists an input sequence $r$
such that
\begin{equation}
\left\vert q_{\Psi}\left[  n\right]  \right\vert \geq |CB|\delta/2,\qquad\forall
n\in\mathbb{Z}_{+}\backslash \{0\}. \label{xPsiBound}%
\end{equation}
Define function $\rho: \mathbb{R}^m \to \mathbb{Z}$ by %Consider the set $\{(2k+1)\delta/2 : k \in \mathbb{Z} \}$, there exists a function $\rho: \mathbb{R}^m \to \mathbb{Z}$ such that
\begin{equation}
\rho(x) = \min \left \{ \arg \min_{k \in \mathbb{Z}} \left[ \displaystyle \frac{2k+1}{2}\delta - (CB)^{-1}CAx \right]\right \}.
\end{equation}
When $r[n]$ is given by
\begin{equation}
r[n] =  \displaystyle \frac{2 \rho(x[n])+1}{2}\delta -(CB)^{-1}CAx[n], \label{BadR}
\end{equation}
we have $r[n] \in[-1,1]$ (since $\delta \in (0,2]$) and
\begin{equation}
|q_\Psi[n+1]|  =  \left| CB\left(\displaystyle \frac{2 \rho(x[n])+1}{2}\delta -u[n]\right) \right | \ge |CB|\delta/2
\end{equation} 
for all $n \in \mathbb Z_+$, because $u[n] \in k\delta$. Hence

\begin{equation}
\mathcal{J}_{G,\phi}(\widehat{\Psi})\geq\phi\left( |CB| \delta/2\right)  . \label{gamma_minLB}
\end{equation}

Inequalities  \eqref{gamma_minUB} and  \eqref{gamma_minLB} complete the proof.

%From the proof we can see that without loss of generality we can assume that $CB =1$. Indeed, since for every $n \in \mathbb Z_+$, $r[n]$ can guarantee that $|q[n+1]| \ge |CB|\delta/2$, monotonicity of function $\phi$ is sufficient for \eqref{gamma_minUB}. Likewise, for every $n \in \mathbb Z_+$, $u[n]$ can guarantee that $|q[n]| \le |CB|\delta/2$, thus monotonicity of function $\phi$ is sufficient for \eqref{gamma_minLB}. Therefore, the case when $CB \ne 0$ can be converted to $CB=1$, since scaling $G(z)$ only scales the performance.
   %
%%%%%%%%%%%%%%%%%%%%%%%%%%%%%%%%%%%%%%%%%%%%%%%%%%%%%%%%%%%%%%%%%

\nocite{PhDBib:Sasha_FSA} \nocite{PhDBib:Mitra0} \nocite{PhDBib:ROC1996} \nocite{PhDBib:Basar95}
\bibliographystyle{IEEEtran}
\bibliography{acompat,My_PhD_Bibliography}

\end{document}